\documentclass[12pt]{article}
\usepackage[dvips]{graphicx}
\usepackage{amsmath,amsfonts,amssymb,amsthm,color}

\usepackage[numbers]{natbib}

\newtheorem{theorem}{Theorem}
\newtheorem{lemma}{Lemma}

\newtheorem{remark}{Remark}
\newtheorem{proposition}{Proposition}

\newtheorem{algorithm}{Algorithm}

\newtheorem{condition}{Condition}

\usepackage{bm}
\bmdefine{\Bt}{t}
\bmdefine{\BX}{X}
\bmdefine{\BY}{Y}
\bmdefine{\BZ}{Z}
\bmdefine{\BB}{B}
\bmdefine{\BM}{M}
\bmdefine{\BD}{D}
\bmdefine{\Bi}{i}
\bmdefine{\Bj}{j}
\bmdefine{\Bk}{k}
\bmdefine{\Bx}{x}
\bmdefine{\By}{y}
\bmdefine{\Bz}{z}
\bmdefine{\Bv}{v}
\bmdefine{\Bw}{w}
\bmdefine{\Bn}{n}
\bmdefine{\Ba}{a}
\bmdefine{\Bb}{b}
\bmdefine{\Bc}{c}
\bmdefine{\Be}{e}
\bmdefine{\Bu}{u}
\bmdefine{\Bp}{p}
\bmdefine{\Bzero}{0}
\bmdefine{\Bone}{1}

\newcommand{\supp}{\mathop{\mathrm{supp}}}

\newcommand{\figcaption}[1]{\def\cptype{Figure} \caption{#1}}

\newcommand{\figref}[1]{Figure \ref{#1}}

\newcommand{\eq}[1]{(\ref{#1})}

\newcommand{\dg}{\text{deg}_{G_w}}
\newcommand{\dw}{\text{deg}_W}

\title{Graver basis for an undirected graph and its application to testing the beta model of random graphs}
\author{Mitsunori Ogawa%
\thanks{Department of 
Mathematical Informatics, 
Graduate School of Information Science and Technology, 
University of Tokyo.},
Hisayuki Hara%
\thanks{
Faculty of Economics, Niigata University}
and Akimichi Takemura\footnotemark[1]\ \thanks{JST CREST}}
\date{November 2011}

\begin{document}
\maketitle

\begin{abstract}
In this paper we give an explicit and algorithmic description of
Graver basis for the toric ideal associated with a simple undirected
graph and apply the  basis for testing the beta model of random graphs by Markov chain Monte Carlo method.
\end{abstract}


\section{Introduction}
\label{sec:introduction}
Random graphs and their applications to the statistical
modeling of complex networks  have been attracting 
much interest in many fields, including statistical mechanics, ecology,
biology and sociology (e.g.\ \citet{Newman-2003},
\citet{Goldenberg-2009}). 
Statistical models for random graphs have been studied since 
\citet{Solomonoff-Rapoport-1951} and \citet{Erdos-Renyi-1960} introduced  
the Bernoulli random graph model.
The beta model generalizes the Bernoulli model to 
a discrete exponential family with 
vertex degrees as sufficient statistics.
The beta model was discussed by \citet{holland-leinhardt-jasa1981} in
the directed case and by \citet{Park-Newman-2004}, 
\citet{blitzstein-diaconis} and \citet{chatterjee-diaconis-sly} in the
undirected case. 
The Rasch model \citep{Rasch1980}, which is a standard model in the item
response theory, is also interpreted as a beta model for undirected
complete bipartite graphs.   
In this article we discuss the random sampling of graphs 
from the conditional distribution in the beta model when 
the vertex degrees are fixed.

In the context of social network the vertices of the graph represent
individuals and their edges represent relationships between individuals.  
In the undirected case the graphs are sometimes restricted to be simple,
i.e., no loops or multiple edges exist. 
The sample size for such cases is at most the number of edges of the
graph and is often small.  
The goodness of fit of the model is usually assessed by large
sample approximation of the distribution of a test statistic.    
When the sample size is not large enough, however, it is desirable to use
a conditional test based on the exact distribution of a test statistic.
For the general background on conditional tests and Markov bases,
see \citet{drton-sturmfels-sullivant}.

Random sampling of graphs with a given vertex degree sequence enables 
us to numerically evaluate the exact distribution of a test statistics for the beta
model.
\citet{blitzstein-diaconis} developed a sequential importance sampling
algorithm for simple graphs which generates graphs through operations 
on vertex degree sequence.
In this article we construct a Markov chain Monte Carlo algorithm for sampling
graphs by using the Graver basis for the toric ideal arising from the
underlying graph of the beta model. 

A Markov basis \citep{diaconis-sturmfels} is often used for sampling from
discrete exponential families. 
Algebraically a Markov basis for the underlying graph of the 
beta model is defined as a set of
generators of the toric ideal arising from the underlying graph of the beta model.
A set of graphs with a given vertex degree sequence is called a fiber
for the underlying graph of the beta model. 
A Markov basis for the underlying graph of the 
beta model is also considered as a set of Markov
transition operators connecting all elements of every fiber. 
\citet{petrovic-etal-Urbana2010} discussed some properties of the toric
ideal arising from the model of \cite{holland-leinhardt-jasa1981} and
provided Markov bases of the model for small directed graphs.
Properties of toric ideals arising from a graph
have been studied in a series of papers by Ohsugi and Hibi
(\cite{ohsugi-hibi-1999aam}, \cite{ohsugi-hibi-1999ja}, 
\cite{ohsugi-hibi-2005jaa}).  

The Graver basis is the set of primitive binomials of the toric ideal.  
Applications of the Graver basis to integer programming
are discussed in \citet{onn-book}.
Since the Graver basis is a superset of any minimal Markov basis, the Graver
basis is also a Markov basis and therefore connects every fiber.  
When the graph is restricted to be simple, however, a Markov basis does
not necessarily connect all elements of every fiber. 
A recent result by \citet{hara-takemura-Urbana2010} implies that the set of
square-free elements of the Graver basis connects all elements of every 
fiber of simple graphs with
a given vertex  degree sequence.
Thus if we have the Graver basis, 
we can sample graphs from any fiber, with or without the restriction that
graphs are simple, in such a way that every graph in the fiber
is generated with positive probability. 

In the sequential importance sampling algorithm of
\cite{blitzstein-diaconis} the underlying graph for the model was 
assumed to be complete, i.e., all the edges have positive probability. 
In our approach we can allow that some edges are absent from the
beginning 
(structural zero edges in the terminology of contingency table
analysis), 
such as the bipartite graph for the case of the Rasch model.  
In fact the Graver basis for an arbitrary graph is obtained by restriction
of the Graver basis for the complete graph 
to the existing edges of $G$ (cf. Proposition 4.13 of \citet{sturmfels1996}).
Moreover our algorithm can be applied not only for sampling simple graphs but also for sampling general undirected graphs without substantial adjustment. 
These are the advantages of the Graver basis.

The Graver basis for small graphs can be computed by a computer algebra
system such as 4ti2 (\citet{4ti2}). 
For even moderate-sized graphs, however, it is difficult to compute the Graver basis via 4ti2 in a practical amount of time.   
In this article we first provide a complete description of the Graver basis for an undirected graph. 
In general the number of elements of the Graver basis is too large.
So we construct an adaptive algorithm for sampling elements from the Graver basis, which is enough for constructing
a connected Markov chain over any fiber.
The recent paper of \citet{reyes-takakis-thoma}
discusses the Graver basis for an undirected graph and 
gives a characterization of the Graver basis.
We give a new description of the Graver basis, which is more
suitable for sampling elements from the Graver basis.

The organization of this paper is as follows. 
In Section \ref{sec:beta-model} we give a brief review on some
statistical models for random graphs and clarify the connection between
the models and toric ideals arising from graphs.
In Section \ref{sec:graver} we provide an explicit description of the 
Graver basis for the toric ideal associated with an undirected
graph. 
Section \ref{sec:algorithm} gives an algorithm for random sampling of
square-free elements of the Graver basis.
In Section \ref{sec:examples} we apply the proposed algorithm to some 
data sets and confirm that it works well in practice.
We conclude the paper with some remarks
in Section \ref{sec:remarks}.

\section{The beta model of random graphs}
\label{sec:beta-model}

In this section we give a brief review of the beta model for undirected graphs
according  to \citet{chatterjee-diaconis-sly}. 

Let $G$ be an undirected graph with $n$ vertices 
$V(G) = \{ 1, 2, \ldots ,n \}$. 
Here we assume that $G$ has no loop.
Let $E=E(G)$ be the set of edges.
For each edge $\{i,j\} \in E$, let a non-negative integer 
$x_{ij}$ be the weight for $\{i,j\}$
and denote $\bm{x} = \{x_{ij} \mid \{i,j\} \in E\}$. 
$\bm{x}$ is considered as an $\vert E \vert$ dimensional integer vector.
We assume that an observed graph $H$ is generated by independent binomial distribution  
$B(n_{ij},p_{ij})$ for each edge $\{i,j\} \in E$, i.e., 
$x_{ij} \sim B(n_{ij},p_{ij})$ with    
\begin{align*}
 p_{ij} & := \frac{e^{\beta_i + \beta_j}}{1 + e^{\beta_i + \beta_j}} 
 = \frac{\alpha_i \alpha_j}{1+ \alpha_i \alpha_j},  \qquad \alpha_i = e^{\beta_i}.
\end{align*}
Then the probability of $H$ is described as
\begin{align}
 \label{beta-model}
 P(H) & \propto \prod_{\{i,j\} \in E}
 p_{ij}^{x_{ij}}(1-p_{ij})^{n_{ij}-x_{ij}}\nonumber\\ 
 & = \frac{1}{\prod_{\{i,j\} \in E} ( 1 + \alpha_i\alpha_j )^{n_{ij}}}
 \prod_{\{i,j\} \in E}(\alpha_i \alpha_j)^{x_{ij}}\nonumber\\
 & = \frac{\prod_{i \in V} \alpha_i^{\sum _{ j: \{ i,j \} \in E } x_{ij} }}{\prod_{\{i,j\} \in E} (1+\alpha_i\alpha_j)^{n_{ij}}}.
\end{align}
The model (\ref{beta-model}) is called the beta model
\citep{chatterjee-diaconis-sly}. 
Note that if $x_{ij}=0$ then the observed graph $H$ does not have an
edge $\{i,j\}$ even if $\{i,j\}\in E(G)$ for the underlying graph $G$.

This model was considered by many authors (e.g. \citet{Park-Newman-2004}, 
\citet{blitzstein-diaconis} and \citet{chatterjee-diaconis-sly}).
The $p_1$ model for random directed graphs by
\citet{holland-leinhardt-jasa1981} can be interpreted as a
generalization of the beta model. 
When $G$ is a complete bipartite graph, the beta model coincides with the Rasch
model \citep{Rasch1980}.  
The many-facet Rasch model by \citet{linacre1989}, which is a
multivariate version of the Rasch model, can be interpreted as a
generalization of the beta model such that $G$ is a complete $k$-partite
graph.  

Let $d_1, \ldots, d_n$ be a degree sequence, i.e., $d_i := \sum_{j: \{i,j\}\in E} x_{ij}$ for each vertex $i$.
Denote $\bm{d}:= (d_1,\ldots,d_n)$.
The sufficient statistic for (\ref{beta-model}) is $\bm{d}$. 
Let $\bm{A}: \vert V \vert \times  \vert E \vert$ denote the
incidence matrix between vertices and edges of $G$.  Then 
it is easily seen that $\bm{x}$ and $\bm{d}$ are related as
\[
 \bm{A}\bm{x}=\bm{d}.
\]

A set of graphs (without restriction to be simple)
${\cal F}_{\bm{d}}=\{ \bm{x}\ge 0 \mid \bm{A}\bm{x}=\bm{d}\}$
with a given degree sequence $\bm{d}$ is called a fiber 
for  $\bm A$ (or for the underlying graph $G$).
An integer array $\bm{z}$ of the same dimension as $\bm{x}$ is called a
move 
if $\bm{A}\bm{z}=0$. 
A move $\bm{z}$ is written as the difference of its positive part and
negative part as $\bm{z} = \bm{z}^+ - \bm{z}^-$.
Since $\bm{A}\bm{z} = \bm{A}\bm{z}^+ - \bm{A}\bm{z}^-$, 
every move is written as the difference of two graphs in the same fiber.
A finite set of moves is called a Markov basis for the incidence matrix $\bm{A}$ 
if for every fiber any two graphs are mutually accessible 
by the moves in the set \citep{diaconis-sturmfels}.
By adding or subtracting moves in a Markov basis, 
we can sample graphs from any fiber in such a way that every graph in
the fiber is generated with positive probability. 
Note that $x_{ij}$ in the beta model
 (\ref{beta-model}) is restricted as $0\le x_{ij}\le n_{ij}$.
We denote the subset of the fiber ${\cal F}_{\bm{d}}$ with this restriction as
${\cal F}_{\bm{d},\bm{n}}=\{\bm{x} \mid \bm{A}\bm{x}=\bm{d}, 0\le x_{ij}\le n_{ij}, \{i,j\}\in E\}$.

\begin{figure}[!h]
\vspace{-3mm}
\begin{center}
\includegraphics[]{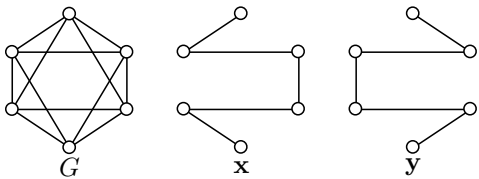}
\vspace{-1mm}
\figcaption{Example graphs.}
\label{fig:example_graph}
\end{center}
\end{figure}

To assess the goodness of fit of the beta model
we usually utilize a large sample approximation of the distribution of a test statistics.
However,
when $n_{ij}$'s are not large enough,
it is not appropriate to use the large sample approximation.
Especially, as mentioned in Section \ref{sec:introduction}, 
graphs are restricted to be simple ($n_{ij}\equiv 1$) in some practical problems. 
For a simple graph, $x_{ij}$, $\{i,j\}\in E$,  is either zero or one.
A Markov basis for the incidence matrix $\bm{A}$
guarantees the connectivity of every fiber ${\cal F}_{\bm{d}}$ if the
restriction that graphs are simple is not imposed. 
Under the restriction, however, a Markov basis does not necessarily
connect the subset ${\cal F}_{\bm{d},\bm{1}}$ of the fiber 
${\cal F}_{\bm{d}}$. 
For example, consider the beta model with the underlying graph $G$
in \figref{fig:example_graph} and $n_{ij}=1$ for each edge $\{ i,j \} \in E$.
It can be shown that 
a set of all $4$-cycles in $G$ is a Markov basis for the incidence matrix of $G$.
However ${\bm x}$ and ${\bm y}$ in  \figref{fig:example_graph} are not mutually accessible by 4-cycles under the restriction that graphs are simple.

For a given $\Bx$, $\supp(\Bx)=\{ e \mid x_e > 0 \}$ denotes the set of
observed  edges of $\Bx$. 
For two moves $\Bz_1, \Bz_2$, the sum $\Bz_1 + \Bz_2$ is called conformal
if there is no cancellation of signs in $\Bz_1 + \Bz_2$, i.e., 
$\emptyset = \supp(\Bz_1^+) \cap \supp(\Bz_2^-)= \supp(\Bz_1^-)\cap \supp(\Bz_2^+)$.
The set of moves which can not be written as a conformal sum of two nonzero moves
is called the Graver basis. 
The Graver basis is known to be a Markov basis 
\citep[e.g.][]{drton-sturmfels-sullivant}.
A move is square-free if the absolute values of its elements are $0$ or $1$.
By the same augment of Proposition 2.1 of \citet{hara-takemura-Urbana2010}, 
we can obtain the following proposition.
\begin{proposition}
 \label{prop:graver}
The Graver basis for the underlying graph of the beta model connects all elements of every fiber.
Furthermore, the set of square-free moves of the Graver basis connects 
all elements of every fiber with the restriction of simple graphs.
\end{proposition}
\begin{proof}
Let $\bm{x}, \bm{y}$ be two elements of the same fiber.
The difference $\bm{y} - \bm{x}$ is written as a conformal sum of primitive moves:
\begin{eqnarray} \label{eqn:conf}
  \bm{y} - \bm{x} = \bm{z}_1 + \cdots + \bm{z}_r
\end{eqnarray}
where $\bm{z}_i, 1 \leq i \leq r, $ are elements of the Graver basis.
Since there is no cancellation of signs on the right hand side,
$\bm{x} + \bm{z}_1 + \cdots + \bm{z}_k$ belongs to the same fiber for 
$k=1,\ldots ,r$.
Therefore the Graver basis connects all elements of every fiber.

Suppose $n_{ij}=1$ for every $ \{ i,j \} \in E$ in the setting of the beta model.
It is easy to see that each $\bm{z}_i$ is square-free in \eq{eqn:conf}.
It means that the set of square-free moves of the Graver basis connects 
all elements of every fiber with the restriction of simple graphs.
\end{proof}

Therefore it suffices to have the Graver basis to sample graphs from any
fiber with or without the restriction that graphs are simple. 
In the next section we derive the Graver basis 
for the underlying graph of the beta model 

\section{Graver basis for an undirected graph}
\label{sec:graver}


In this section we will give two characterizations of the Graver basis for an undirected graph.
Theorem \ref{theorem:graver2} in Section \ref{subsec:characterization}
is the main result of this paper which gives a necessary and sufficient condition 
for a element of the Graver basis as a sequence of vertices.
Proposition \ref{proposition:graver1}, which is used 
for the proof of Theorem \ref{theorem:graver2},
gives a characterization of the Graver basis through recursive operations on
the graph, which is of some independent interests.

\subsection{Preliminaries}
\label{subsec:definitions}
Let $G=(V(G),E(G))$ be a simple connected graph
with $V(G)= \{ 1, 2, \ldots ,n \}$
and $E(G) = \{ e_1, e_2,\ldots, e_m \}$.
A {\it walk} connecting $i \in V(G)$ and $j \in V(G)$ is a finite sequence of edges of the form
\begin{eqnarray*}
	w = ( \{ i_1,i_2 \}, \{ i_2,i_3 \}, \ldots , \{ i_q, i_{q+1} \} )
\end{eqnarray*}
with $i_1=i,i_{q+1}=j$. 
The {\it length} of the walk $w$ is the number  of edges $q$ of the walk.
An {\it even} (respectively {\it odd}) walk is a walk of even (respectively {\it odd}) length.
A walk $w$ is {\it closed} if $i=j$.
A {\it cycle} is a closed walk
$w = ( \{ i_1,i_2 \}, \{ i_2,i_3 \}, \ldots , \{ i_q, i_1 \} )$
with $i_l \neq i_{l'}$ for every $1 \leq l < l' \leq q$.

For a walk $w$, let $V(w)=\{i_1, \dots, i_{q+1}\}$ 
denote the set of vertices appearing in $w$ and 
let $E(w)=\{\{ i_1,i_2 \}, \{ i_2,i_3 \}, \ldots , \{ i_q, i_{q+1} \}\}$
denote the set of edges appearing in $w$.
Furthermore let $G_w=(V(w),E(w))$
be the subgraph of $G$, whose vertices and edges appear in the walk $w$. 

In order to describe known results on the toric ideal $I_G$
arising from an undirected graph $G$, 
we give an algebraic definition of $I_G$.
Let $K[{\bf t}]=K[t_1, \ldots, t_n]$ be a polynomial ring
in $n$ variables over $K$.
We will associate each edge $e_r=\{ i , j \} \in E(G)$
with the monomial ${\bf t}_r=t_it_j \in K[{\bf t}]$.
Let $K[{\bf s}]=K[s_1, \ldots, s_m]$ be a polynomial ring
in $m=|E(G)|$ variables over $K$ 
and let $\pi$ be a homomorphism from $K[{\bf s}]$ to $K[{\bf t}]$ 
defined by $\pi: s_r \mapsto {\bf t}_r$.
Then the {\it toric ideal} $I_G$ of the graph $G$ is defined as
\begin{eqnarray*}
	I_G = \ker (\pi) = \{ f \in K[{\bf s}] \mid \pi (f) = 0 \}.
\end{eqnarray*}
A binomial $f = u-v\in I_G$ is called {\it primitive} if there is no binomial $g = u^{\prime}-v^{\prime}\in I_G$, $g \neq 0 ,f$, such that $u^{\prime} | u$ and $v^{\prime} | v$.
The {\it Graver basis} of $I_G$ is the 
set of all primitive binomials belonging to $I_G$ and we denote it by ${\cal G} (I_G)$.  
If we write  the monomials $u,v$ as $u=s^{\bm{x}}, v=s^{\bm{y}}$,
then $u-v\in I_G$ if and only if $\bm{x} - \bm{y}$ is a move. Furthermore
$u-v\in I_G$ is primitive if and only if $\supp(\bm{x})\cap \supp(\bm{y})=\emptyset$ and
$\bm{x} - \bm{y}$ can not be written as a conformal sum of two nonzero moves.

For a given even closed walk
$w = (e_{j_1}, e_{j_2}, \ldots , e_{j_{2p}})$
we define a binomial $f_w \in I_G$ as
\[
f_w = f_w^+ - f_w^- , \qquad  \text{where}\ \ 
f_w^+ = \prod _{k=1}^p s_{j_{2k-1}},~~
f_w^- = \prod _{k=1}^p s_{j_{2k}}.
\]
An even closed walk $w^{\prime}$ is a {\it proper subwalk} of $w$, if $g_{w^{\prime}}^+ \mid f_w^+$ and $g_{w^{\prime}} ^- \mid f_w^-$ hold for the binomial $g=g_{w^{\prime}}^+ - g_{w^{\prime}}^- (\neq f_w)$.
Note that even if there is a proper subwalk $w^{\prime}$ 
of an even closed walk $w$,
$w^{\prime}$ dose not necessarily go along with $w$,
i.e., the edges of $w^{\prime}$ may not appear as consecutive edges of $w$.
An even closed walk $w$ is called {\it primitive},
if its binomial $f_w$ is primitive.
Then the primitiveness of $w$ is equal to non-existence of a proper subwalk of $w$.

A  characterization of the primitive walks of graph $G$, which gives a 
necessary condition for a binomial to be primitive, 
was given by \citet[]{ohsugi-hibi-1999ja}.
\begin{proposition}[\cite{ohsugi-hibi-1999ja}]
\label{prop:hibi1}
Let $G$ be a finite connected graph.
If $f \in I_G$ is primitive, then we have $f=f_w$ where $w$ is one of the following even closed walks:
\begin{enumerate}
\renewcommand{\labelenumi}{(\roman{enumi})}
\setlength{\itemsep}{0pt}
\item $w$ is an even cycle of $G$.
\item $w = (c_1,c_2)$, where $c_1$ and $c_2$ are odd cycles of $G$ having exactly one common vertex.
\item $w = (c_1,w_1,c_2,w_2)$, where $c_1$ and $c_2$ are odd cycles of $G$ having no common vertex and where $w_1$ and $w_2$ are walks of $G$ both of which contain a vertex $v_1$ of $c_1$ and a vertex $v_2$ of $c_2$.
\end{enumerate}
\end{proposition}
Every binomial in the first two cases is primitive but a binomial in the third case is
not necessarily primitive.

\subsection{Characterization of primitive walks}
\label{subsec:characterization}

\begin{figure}[!h]
\begin{center}
\includegraphics[]{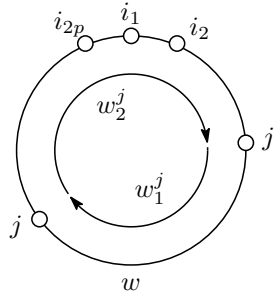}
\vspace{-1mm}
\figcaption{Even closed walk $w$.}
\label{fig:condition_graver2}
\end{center}
\end{figure}

In this subsection we give a simple characterization of the primitive walks of a graph $G$ as sequences of vertices.
Express an even closed walk $w$ as a sequence of vertices: $(i_1, i_2 , \ldots ,i_{2p}, i_1)$,
where $i_1 \equiv i_{2p+1}$.
Let $\#_w(i)=\#\{1 \le l \le 2p \mid i_l= i\}$ 
denote the number of times $i$ is visited in the walk $w$ before it returns to the vertex $i_1$. 
Consider the following condition for the even closed walk $w$.
\begin{condition}\label{condition:graver2}
(i) $\#_w(i) \in \{ 1,2 \}$ for every vertex $i \in V(w)$.
(ii) For every vertex $j \in V(w)$ with $\#_w(j) = 2$ and 
$j=i_l = i_{l'}$, $1 \le l < l' \le 2p$, 
the closed walks $w_1^j = (i_l, \ldots ,i_{l'})$ and $w_2^j = (i_{l'}, \ldots, i_{2p}, i_1, \ldots ,i_{l-1}, i_l)$ are odd walks with $V(w_1^j) \cap V(w_2^j) =\{ j \}$.
(cf.\ \figref{fig:condition_graver2}).
\end{condition}
\begin{remark}
The equation $V(w_1^j) \cap V(w_2^j) =\{ j \}$ in Condition \ref{condition:graver2} means that there are no crossing chords in \figref{fig:condition_graver2} when adding a chord $\{ j, j \}$ to the figure for every vertex $j \in V(w)$ with $\#_w(j) = 2$.
\end{remark}

Using Condition \ref{condition:graver2}, we can characterize the Graver basis for a graph $G$ as follows.
\begin{theorem}\label{theorem:graver2}
A binomial $f \in I_G$ is primitive if and only if there exists an even closed walk $w$ with $f_w=f$ satisfying Condition \ref{condition:graver2}. 
\end{theorem}
\begin{remark}
It follows from the definition of primitive walks and Theorem \ref{theorem:graver2} that
if an even closed walk $w$ is primitive,
every even closed walk $w^{\prime}$ with $f_{w^{\prime}} = f_w$ is primitive and satisfies Condition \ref{condition:graver2}. 
\end{remark}
\begin{remark}
As mentioned in Section \ref{sec:introduction},
there is another characterization of Graver basis in Theorem 3.1 of
\citet{reyes-takakis-thoma}.
It also gives a necessary and sufficient condition for the 
primitiveness of even closed walks,
by using some new graphical concepts such as {\it ``block''} and {\it ``sink''}. 
Our characterization in Theorem \ref{theorem:graver2} gives a simpler description of Graver basis, because it does not need any new graphical concepts.
Furthermore it is more convenient in the algorithmic viewpoint: 
When an even closed $w$ is given as a sequence of vertices or edges, 
we can easily determine if $w$ is primitive by checking directly
Condition \ref{condition:graver2} 
without distinguishing any graphical objects.
\end{remark}

Before proving Theorem \ref{theorem:graver2}, 
we state another characterization of primitive walks given in Proposition
\ref{proposition:graver1} below.
In order to that, we need some more definitions on graphs.
For a walk $w = (e_{j_1}, e_{j_2}, \ldots , e_{j_{q}})$, let $W=W(w)$ denote the weighted subgraph $(V(w),E(w),\rho)$ of $G$ where $\rho : E(w) \rightarrow {\mathbb Z}$ is the weight function defined by $\rho (e) := \#\{ l \mid e_{j_{2l+1}}= e\} - \#\{ l \mid e_{j_{2l}}= e\}$ for each edge $e \in E(w)$.
For simplicity, we denote a weight $+1$ (respectively $-1$) by $+$ (respectively $-$) in our figures.
For a vertex $i \in V(w)$, we define two kinds of degrees of vertex $i$:
\begin{eqnarray*}
	\dg (i) &=& \#\{ e \in E(w) \mid i \in e \} , \\
	\dw (i) &=& \sum_{e \in E(w) : i \in e} |\rho (e)|.
\end{eqnarray*}
$\dg (i)$ is the usual degree of $i$ in $G_w$.  Note that the same weighted graph
$W$ might correspond to two different even closed walks $w,w'$, i.e.\ 
$W(w)=W(w')$. Given a weighted graph $W$, we say that $w$ {\it spans} $W$
if $W=W(w)$ and $\{ e_{j_l} \mid l\text{:odd} \} \cap \{ e_{j_l} \mid l\text{:even} \} = \emptyset $.

\begin{figure}[!h]
\begin{center}
\includegraphics[]{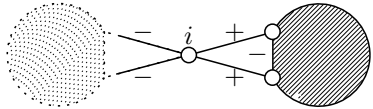}
\vspace{-0mm}
\figcaption{A vertex $i$ whose separation is not defined.}
\label{fig:nonseparation}
\end{center}
\end{figure}

Now we define two operations, {\it contraction} and {\it separation}, 
on a weighted graph $W$.
\begin{itemize}
 \item  Let $e=\{ i,j \} \in E(w)$ be an edge with $|\rho (e)|=2$, 
whose removal from $G_w$ increases the number of connected components of the remaining subgraph.
			{\it Contraction} of $e$ is an operation as shown in \figref{fig:reduction}.
			That is, it first replaces $W$ by $W^{\prime } = (V(w) \setminus \{ j \}, E^{\prime}, \rho ^{\prime})$ where $E^{\prime}$ consists of all edges of $W$ contained in $V(w) \setminus \{ j \}$, together with all edges $\{ \alpha , i \}$, where $\{ \alpha , j \}$ is an edge of $W$ different from $e$.
			Then, it defines $\rho ^{\prime}$ by inversion of the signs of weights of edges belonging  to the $i$-side of $W$.
\begin{figure}[!h]
\begin{center}
\includegraphics[width=0.7\textwidth ]{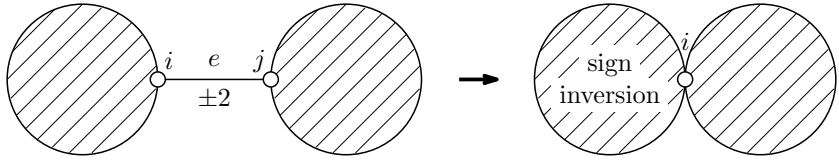}
\vspace{-0mm}
\figcaption{Contraction.}
\label{fig:reduction}
\end{center}
\end{figure}
 \item  Let $i \in V(w)$ be a vertex with $\dg(i)=\dw(i)=4$, such that the removal of $i$ increases the number of connected components of the remaining subgraph and the positive side as well as the negative side of $i$ fit to one of three cases (a)--(c) (respectively to the sign reverse cases) in \figref{fig:separation1}.
			{\it Separation} of $i$ is an operation as shown in \figref{fig:separation1}.
			That is, it first deletes the vertex $i$ and all edges connected to $i$ on $W$.
			Then, in the case of (a),  it adds a new edge $\{ k_1, k_2 \}$ with weight $+1$.
			In the case of (b), it redefines $\rho (\{ k_1, k_2 \}) :=+2$ and then contracts $\{ k_1, k_2 \}$, where we assume that the contraction
of $\{k_1, k_2\}$ is possible.
			In the case of (c), it redefines $\rho (\{ k_1, k_2 \}) :=0$.
                        We call this $\{k_1, k_2\}$ an edge with weight $0$.
			The sign reverse cases are defined in the same way.
\begin{figure}[!h]
\begin{center}
\includegraphics[width=0.85\textwidth]{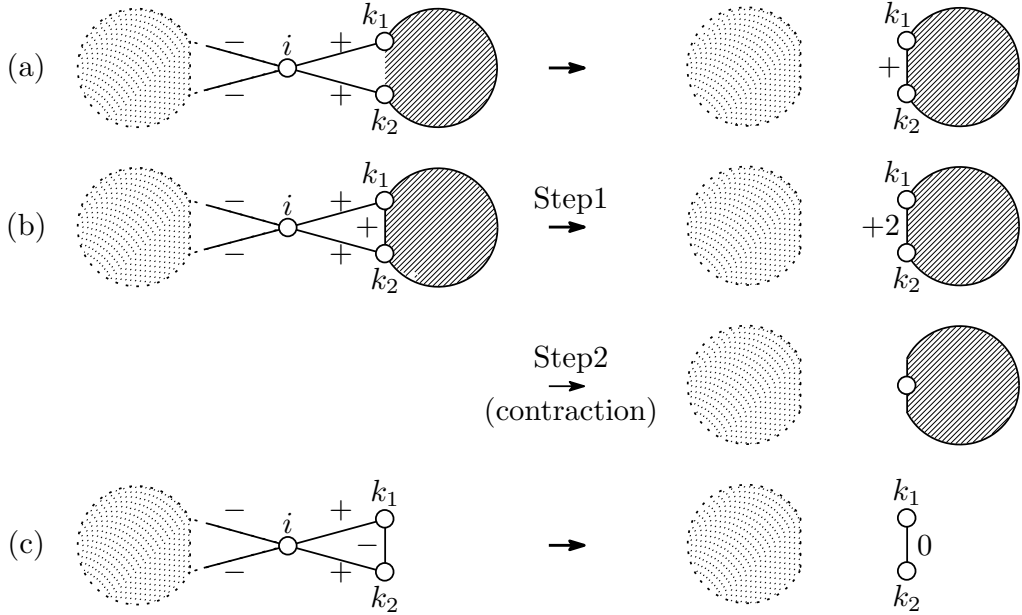}
\vspace{-0mm}
\figcaption{Separation.}
\label{fig:separation1}
\end{center}
\end{figure}
\end{itemize}
Note that the separation is not defined for any vertex $i$ with $\dg (i) = \dw (i) =4$,
if $i$ fits to none of three cases (a)--(c) in Figure \ref{fig:separation1}.
The vertex $i$ in \figref{fig:nonseparation} is such an example, because its positive side fits to none of three cases (a)--(c) in \figref{fig:separation1}.

Let {\it insertion} and {\it binding} be the reverse operations 
of contraction and separation, respectively.
With these operations, consider the following condition for an even closed walk $w=(e_{j_1},e_{j_2},\ldots ,e_{j_{2p}})$.
\begin{condition}\label{condition:graver}
(i) $\{ e_{j_l} \mid l\text{:odd} \} \cap \{ e_{j_l} \mid l\text{:even} \} = \emptyset $.
Every vertex $i \in V(w)$ satisfies $\dw (i)\in \{ 2,4 \}$.
For every vertex $i$ with $\dw (i)=4$, its removal from $G_w$ increases the number of connected components of the remaining subgraph.
(ii)
Let $\tilde{W}$ be a graph obtained by recursively applying contraction and separation of all possible edges and vertices in $W$.
Then each connected component of $\tilde{W}$ is an even cycle or an edge with weight 0.
\end{condition}
\begin{proposition}\label{proposition:graver1}
For an even closed walk $w$,
the binomial $f_w$ is primitive 
if and only if $w$ satisfies Condition \ref{condition:graver}.
\end{proposition}

We establish some lemmas to prove Proposition \ref{proposition:graver1}.
Our proof also shows that $\tilde{W}$ in Condition \ref{condition:graver}
does not depend on the order of application of contractions and separations excepting the sign inversion of weights of edges of each connected component in $\tilde{W}$.
\begin{lemma}\label{lemma:nc1}
If $w$ is a primitive walk,
$w$ satisfies (i) in Condition \ref{condition:graver}.
\end{lemma}
\begin{proof}
Consider a vertex $i \in V(w)$.
Since $w$ is closed, $\dw (i)$ is even.
Furthermore, since $w$ is primitive, $\{ e_{j_l} \mid l\text{:odd} \} \cap \{ e_{j_l} \mid l\text{:even} \} = \emptyset $ holds
which implies that there is no cancellation in the calculation of weight on any edge.
Then, a half of the weight $\dw (i)/2$ is  assigned as positive weights and other half $\dw (i)/2$  is assigned as negative weights to the edges connected  to $i$ on $W$.
Therefore $\dw (i) \in \{ 2,4,6,\ldots \}$.
Now suppose $\dw (i) \geq 6$.
Consider that we start from a vertex $i$ along an edge with positive weight and go along the walk $w$ or its reverse until returning back to $i$ again for the first time. 
Since $w$ is primitive, we have to come back to $i$ along an edge with positive weight for the first time.
Let us continue along $w$ or its reverse until returning back to $i$.
By the same reasoning, the last edge of this closed walk has a negative weight.
This implies that this even closed walk becomes a proper subwalk of $w$, a contradiction to the primitiveness of $w$.
Therefore $\dw (i)$ is $2$ or $4$.

To prove the remaining part, let $i \in V(w)$ be a vertex with $\dw (i) =4$ and
consider all closed walks on $W$, where the edge starting from $i$ and
the edge coming back to $i$ have positive weights.
Let $V^+$ be the set of vertices other than $i$ which appear in one of these walks and $V^-$ is defined in the same way.
Then $V^+\cup V^- \cup \{ i \} = V(w)$ holds.
First, we show $V^+ \cap V^- = \emptyset$.
Suppose that there exists a vertex $j \in V^+ \cap V^-$.
Then, as shown in \figref{fig:inW1}, 
there are two closed walks $( \{ i,i_1^+ \}, \Gamma _1^+, \Gamma _2^+, \{ i_2^+,i \} )$ and $( \{ i,i_1^- \}, \Gamma _1^-, \Gamma _2^-, \{ i_2^-,i \} )$.
This implies that we can construct a proper subwalk of $w$ by the combination of $\{ i,i_k^+ \}, \Gamma_k^+(k=1,2)$, and $\Gamma_l^-, \{ i_l^-,i \} (l=1,2)$, a contradiction to the primitiveness of $w$.
Therefore $V^+ \cap V^- = \emptyset$.
\begin{figure}[!h]
\begin{center}
\includegraphics[]{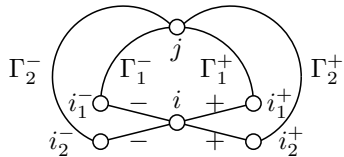}
\vspace{-1mm}
\figcaption{Case that there exists a vertex $j \in V^+ \cap V^-$.}
\label{fig:inW1}
\end{center}
\end{figure}
Second, suppose that the removal of the vertex $i$ from $G_w$ does not increase the number of connected components of the remaining subgraph.
Then, there are vertices $v^+ \in V^+, v^- \in V^-$ such that $\{ v^+,v^- \} \in E(w)$, because $V^+ \cap V^- = \emptyset$ holds as shown above.
Hence, as shown in \figref{fig:inW2}, an even closed walk $( \{ i,i_k^+ \} , \Gamma_k^+, \{ v^+,v^- \} , \Gamma_l^-, \{ i_l^-,i \})$ is a proper subwalk of $w$ for appropriate $k,l \in \{1,2\}$, $k\neq l$, which contradicts to the primitiveness of $w$.
Therefore the removal of $i$ from $G_w$ increases the number of connected components of the remaining subgraph.
\begin{figure}[!h]
\begin{center}
\includegraphics[]{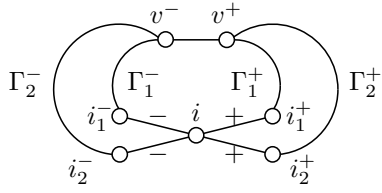}
\vspace{-1mm}
\figcaption{Case that there exists an edge $\{ v^+,v^- \}$.}
\label{fig:inW2}
\end{center}
\end{figure}
\end{proof}

In the following four lemmas,
we show that contraction, separation, and these inverse operations preserve 
the primitiveness of an even closed walk.
The proofs of lemmas are postponed to Appendix.
\begin{lemma}\label{lemma:nc2}
Let an even closed walk $w$ be primitive
and $\tilde{W}$ be the weighted graph which is obtained by a contraction for an edge with its weight $\pm 2$ on $W$.
Then any even closed walk $\tilde{w}$ spanning $\tilde{W}$ is primitive.
\end{lemma}
\begin{lemma}\label{lemma:nc3}
Let an even closed walk $w$ be primitive
and $W_1,W_2$ be the weighted graphs obtained by the separation of a vertex $i$.
Then any even closed walks $w_l (l=1,2)$ spanning $W_l (l=1,2)$
are  primitive or of length two with $f_{w_l}=0$.
\end{lemma}
\begin{lemma}\label{lemma:sc2}
Let $w$ be a primitive walk and let $\tilde{W}$ be the weighted graph obtained by the insertion to $i$ with $\dw (i)=4$ on $W$.
Then any even closed walk $\tilde{w}$ spanning $\tilde{W}$ is primitive.
\end{lemma}
\begin{lemma}\label{lemma:sc1}
Let each $w_l ~(l=1,2)$ be a primitive walk or a closed walk with length two,
and $W$ be the weighted subgraph obtained by binding of $W_1$ and $W_2$.
Then any even closed walk $w$ spanning $W$ is primitive.
\end{lemma}

We now give proofs of Proposition \ref{proposition:graver1}
and 
Theorem \ref{theorem:graver2}.

\begin{proof}[Proof of Proposition \ref{proposition:graver1}]
Let $w$ be a primitive walk.
From Lemma \ref{lemma:nc1} $w$ satisfies (i) in Condition \ref{condition:graver} and every edge $e$ with $|\rho (e)|=2$ can be contracted.
Furthermore, it is easy to see that every vertex $i$ with $\dw(i)=4$ can be separated after recursively applying contractions of all possible edges.
Therefore $\dw(i)=2$ holds for every vertex $i$ on $\tilde{W}$.
From Lemmas \ref{lemma:nc2} and \ref{lemma:nc3}, each even closed walk corresponding to the connected component of $\tilde{W}$ is primitive or of length two.
Then, every connected component of $\tilde{W}$ is an even cycle or an edge with weight 0, because from Proposition \ref{prop:hibi1} every primitive walk includes a vertex $i$ with $\dw (i) =4$ if it is not an even cycle.
Therefore, a primitive walk $w$ satisfies Condition \ref{condition:graver}.
Conversely, suppose an even closed walk $w$ satisfies Condition \ref{condition:graver}.
From Proposition \ref{prop:hibi1} and Lemmas \ref{lemma:sc2} and \ref{lemma:sc1},
$w$ is primitive.
\end{proof}

\begin{proof}[Proof of Theorem \ref{theorem:graver2}]
Let $w$ be a primitive walk.
From Lemma \ref{lemma:nc1}, $\#_w(i) \in \{ 1,2 \}$ holds for each vertex $i \in V(w)$ and $V(w_1^j) \cap V(w_2^j) =\{ j \}$ holds for each vertex $j \in V(w)$ with $\#_w(j) = 2$.
By the primitiveness of $w$, 
the closed walks $w_1^j = (j, \ldots ,j)$ and $w_2^j = (j, \ldots , i_1, \ldots ,j)$ along $w$ are odd closed walks.
Therefore $w$ satisfies Condition \ref{condition:graver2}.

Conversely, let $w$ be an even closed walk with Condition \ref{condition:graver2}.
From Proposition \ref{proposition:graver1}, it suffices to show that $w$ satisfies Condition \ref{condition:graver}.
The condition (i) in Condition \ref{condition:graver} follows from Condition \ref{condition:graver2}.
Then, it is enough to confirm that $w$ satisfies the condition (ii) in Condition \ref{condition:graver}.

First, we claim that every edge $e \in E(w)$ with $|\rho (e)|=2$ can be contracted and every vertex $j$ with $\#_w(j) = 2$ and $\dg (j)=4$, i.e. $\dg(j)=\dw(j)=4$, can be separated.
The case of contraction is obvious from Condition  \ref{condition:graver2}.
We confirm the case of separation.
Consider the vertex $j$ in \figref{fig:test_proof71}.
\begin{figure}[!h]
\begin{center}
\includegraphics[]{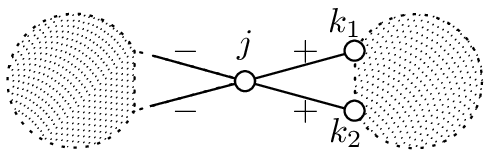}
\vspace{-3mm}
\end{center}
\figcaption{A vertex $j$ with $\dg(j)=\dw(j)=4$.}
\label{fig:test_proof71}
\end{figure}
If an edge $\{ k_1,k_2 \}$ dose not exist or exists with weight $+1$, it belongs to the case (a) or (b) in \figref{fig:separation1}, respectively.
Let us consider the case that there exists an edge $\{ k_1,k_2 \}$ with weight $-1$ and suppose that the vertex $k_1$ connects to more than three edges as shown in \figref{fig:test_proof7}.
\begin{figure}[!h]
\begin{center}
\includegraphics[]{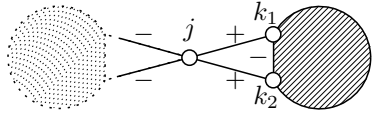}
\vspace{-3mm}
\end{center}
\figcaption{A vertex $j$ which does not exist in $w$ with Condition \ref{condition:graver2}.}
\label{fig:test_proof7}
\end{figure}
Then, $j,k_1$ and $k_2$ appear in $w$ like $(j, k_1, \ldots, k_1,k_2,j)$ or $(j, k_1, \ldots, k_1, k_2,\ldots, k_2,j)$, because $V(w_1^j) \cap V(w_2^j) =\{ j \}$ holds.
This implies that $(k_1, \ldots, k_1)$ is even as shown in \figref{fig:test_proof73}, which contradicts  Condition \ref{condition:graver2}.
Hence the case with $\{ k_1,k_2 \}$ with weight $-1$ belongs to (c) in \figref{fig:separation1}.
Therefore the claim is confirmed.
\begin{figure}[!h]
\begin{center}
\includegraphics[]{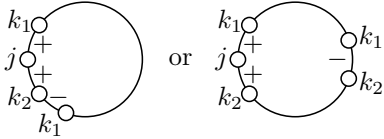}
\vspace{-5mm}
\end{center}
\figcaption{Case that there exists a vertex $j$ in Figure \ref{fig:test_proof7}.}
\label{fig:test_proof73}
\end{figure}

Second, we verify that contraction and separation on $W$ preserve Condition \ref{condition:graver2}.
Consider the case of contraction of an edge $\{i, j\} \in E(W)$.
From Condition \ref{condition:graver2}, such $i,j$ appear in $w$ as $w=(i_1,\ldots,i_{l_1},i,j,i_{l_2},\ldots,i_{l_3}, j,i,i_{l_4},\ldots,i_1)$.
The contraction of $\{i, j\}$ is equivalent to replacing $w$ by $(i_1,\ldots,i_{l_1},i,i_{l_2},\ldots, i_{l_3}, i,$
$i_{l_4},\ldots,i_1)$.
This change causes the decrease of two edges from $w$, and preserves Condition \ref{condition:graver2}.
The case of separation is checked in the same way.

Finally, consider the weighted graph $W^{\prime}$ obtained by all possible contractions and separations on $W$.
From the claims above, every connected component of $W^{\prime}$ satisfies Condition \ref{condition:graver2} and has no vertex $j$ with $\#_w(j) = 2$, i.e.\ an even cycle or an edge with weight 0.
Therefore $w$ satisfies Condition \ref{condition:graver}.
\end{proof}

\section{Algorithm for generating elements of Graver basis}
\label{sec:algorithm}

In this section we present an algorithm for generating elements randomly from the Graver basis for a simple undirected graph.
As shown in Proposition \ref{prop:graver},
for testing the beta model of random graphs with $n_{ij}=1$,
we only need square-free elements of the Graver basis.
Therefore the restriction to square-free elements of our algorithm will be discussed in Remark \ref{remark:algo}.
Theorem \ref{theorem:graver2} guarantees the correctness 
of our algorithm.

We need  some tools in order to construct an algorithm.
Let $T$ be a weighted tree $(V(T),E(T), \mu)$ where $\mu : V(T) \rightarrow {\mathbb  Z}_{\geq 2}=\{2,3,\dots\}$ is a weight function.
For this weighted tree $T$, let us consider the following condition.
\begin{condition}\label{condition:tree}
For each vertex $v_T \in V(T)$, $\deg (v_T) \leq \mu (v_T)$ 
and $\deg (v_T) \equiv  \mu (v_T) \mod 2$.
\end{condition}

With these tools, let us consider generating an element of the Graver basis for a simple undirected graph $G=(V(G),E(G))$.
For simplicity, first consider the case that that $G$ is complete.
We call an edge $e$ with $| \rho (e) | =2$ a cycle in $G_w$ 
for an even closed walk $w$ in this section.
We will discuss later the case that $G$ is not complete.
Let $T=(V(T),E(T),\mu )$ be a weighted tree satisfying Condition \ref{condition:tree} and the following equation:
\begin{eqnarray}\label{eq:vertex_of_T}
	\sum_{v_T \in V(T)}\mu (v_T)-|E(T)|\leq |V(G)|.
\end{eqnarray}
Then, we can construct a primitive walk in $G$ using $T$ as follows.
First, we assign the set of vertices $V_{v_T} \subseteq V(G)$ with $|V_{v_T}| = \mu (v_T)$ for each vertex $v_T \in V(T)$ under the equation
\begin{eqnarray*}
	| V_{v_T} \cap V_{v_T^{\prime}} | = 
	\left\{ \begin{array}{ll}
		1, & \text{if } \{ v_T,v_T^{\prime} \} \in E(T) , \\
		0, & \text{if } \{ v_T,v_T^{\prime} \} \notin E(T), \\
	\end{array} \right.
	\hspace{8mm} (v_T^{\prime} \in V(T))
\end{eqnarray*}
and every vertex $v \in V(G)$ is assigned at most twice.
Equation \eq{eq:vertex_of_T} guarantees that this assignment is possible.
Second, we make cycles in $G$ by arbitrarily ordering the vertices $V_{v_T}$.
Then we make a subgraph of $G$ by taking the union of these cycles.
Finally, we obtain a closed walk by choosing a root vertex from this subgraph and going around it.
It is easy to see that this closed walk is primitive by Theorem \ref{theorem:graver2}.

Conversely we can construct a weighted tree with Condition \ref{condition:tree} and \eq{eq:vertex_of_T} from each primitive walk.
Let $w$ be a primitive walk.
First,  the vertex set $V(T)$ is constructed by creating a vertex $v_c$ of $T$ for each cycle $c$ in $G_w$.
Second, the edge set $E(T)$ is obtained by adding edge $\{ v_c, v_{c^{\prime}} \}$ to $E(T)$ for each pair of cycles $c,c^{\prime}$ in $G_w$ with $V(c) \cap V(c^{\prime}) \neq \emptyset$.
Then, we assign weight $\mu (v_c) := |V(c)|$ to each vertex $v_c \in V(T)$.

Therefore, once we have a weighted tree $T$ with Condition \ref{condition:tree} and \eq{eq:vertex_of_T}, we can construct an element of the Graver basis for $G$.
Such a tree $T$ is constructed by the following algorithm.
\begin{algorithm}[Algorithm for constructing an weighted tree]\ \\
\label{algo:WT}
Input~:~A complete graph $G=(V(G),E(G))$.\\
Output~:~A weighted tree $T=(V(T),E(T),\mu )$ with Condition \ref{condition:tree} and \eq{eq:vertex_of_T}.
\begin{enumerate}
\setlength{\itemsep}{0pt}
	\item  Let $V(T),E(T)$ be empty sets and $n:=|V(G)|$.
	\item  Add a root vertex $r$ to $V(T)$.
	\item  Assign $\mu (r)$ a weight from $\{ 2,3,\ldots,n \}$ randomly.
	\item  Grow $T$ by the following loop.
			  \begin{enumerate}\setlength{\itemsep}{0pt}
			  \renewcommand{\labelenumi}{(\roman{enumi})}
				\item For each vertex $v_T \in V(T)$ which is deepest from $r$, add edges $\{ v_T, v_T^i \}$ to $E(T)$ and the endpoints $v_T^i$ ($i=0,1,\ldots ,I_{v_T}$) to $V(T)$, where the number $I_{v_T}$ is randomly decided under the following two conditions:
						\begin{itemize}\setlength{\itemsep}{0pt}
							\item $I_{v_T} +1 \leq \mu (v_T)$.
							\item $I_{v_T} +1 \equiv \mu (v_T) \mod 2$.
						\end{itemize}
				\item For each new vertex $v_T^i$,  assign $\mu (v_T^i)$ a weight from $\{ 2,3,\ldots,n-\alpha \}$ randomly, where $\alpha := \sum_{v_T \in V(T)} \mu (v_T) -|E(T)|$.
				\item Recompute $\alpha$ and if $\alpha > n$, delete all new vertices and edges in the above (a) and break the loop.
				\item If the total number of new edges is equal to 0, break the loop. 
				\item Return to (a).
			  \end{enumerate}
	\item  If $|V(T)|=1$ and $\mu (r)$ is odd,  change $\mu (r)$ to $\mu (r)-1$ or $\mu (r)+1$.
	\item  If $|V(T)|>1$ and $T$ has a leaf with even weight, subtract or add 1 to the weight.
	\item Output $T$.
\end{enumerate}
\end{algorithm}

Algorithm \ref{algo:WT} provides a simple algorithm for generating an element of Graver basis as follows.
\begin{algorithm}[Algorithm for generating an element of Graver basis]\ \\
\label{algo:TP}
Input~:~A complete graph $G=(V(G),E(G))$.\\
Output~:~A primitive walk $w$.
\begin{enumerate}
\setlength{\itemsep}{0pt}
	\item  Construct a weighted tree $T$ with Condition \ref{condition:tree} and \eq{eq:vertex_of_T} by Algorithm \ref{algo:WT}.
	\item  Construct a primitive walk by assigning vertices of $G$ and ordering them randomly.
	\item Output $w$.
\end{enumerate}
\end{algorithm}

\begin{figure}[!h]
\begin{center}
\includegraphics[width=0.95\textwidth]{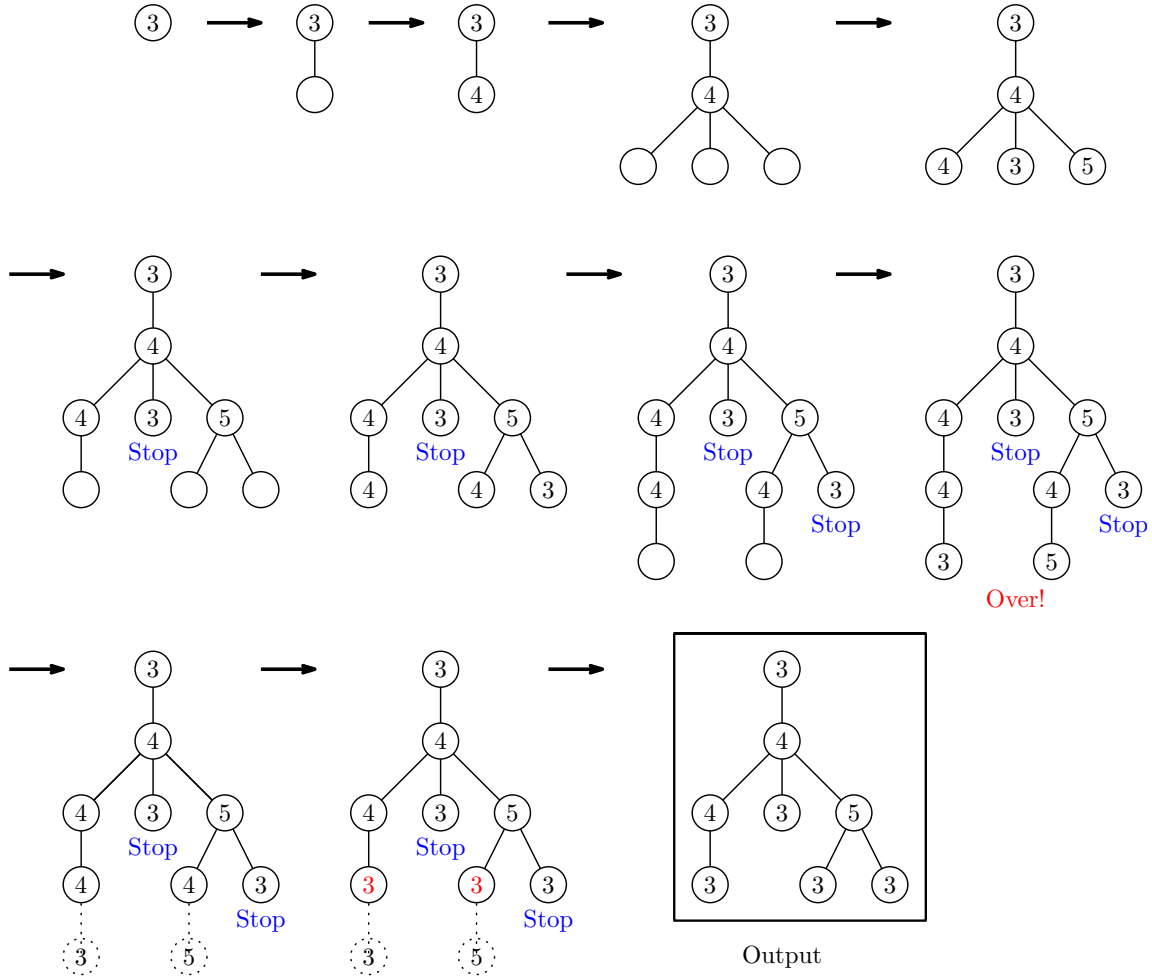}
\vspace{-3mm}
\end{center}
\figcaption{Demonstration of Algorithm \ref{algo:WT}.}
\label{fig:tree_algo}
\end{figure}

\begin{figure}[!h]
\begin{center}
\includegraphics[width=9cm]{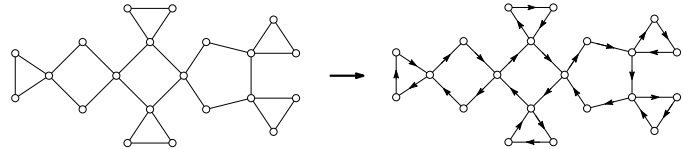}
\vspace{-3mm}
\end{center}
\figcaption{Demonstration of Algorithm \ref{algo:TP}.}
\label{fig:tree_algo_graph}
\end{figure}

Since there is no restarts in Algorithm \ref{algo:TP}, it has a fixed worst case running time for a complete graph $G$.
In each step, the algorithm performs $O(|V(G)|)$ operations.
Then it generates one element of the Graver basis for $G$ in $O(|V(G)|)$ time.
A demonstration for the case of a complete graph $G$ with $|V(G)|=25$ is shown in Figures \ref{fig:tree_algo} and \ref{fig:tree_algo_graph}.
The output of this demonstration is a primitive walk $w$ with $|V(G_w)| = 21$
in Figure \ref{fig:tree_algo_graph}.

\begin{remark}
\label{remark:algo}
For the case that an input graph $G$ is not complete,
the elements of the Graver basis for $G$ can be generated by throwing away elements with supports not contained in $G$
(Proposition 4.13 of \citet{sturmfels1996}). In fact this is the advantage
of considering the Graver basis.
The restriction for generation of square-free elements of the Graver basis 
can be realized by a slight modification in Algorithm \ref{algo:WT}.
In fact, it suffices to change merely $\{ 2,3,4,\ldots \}$ to $\{ 3,4,\ldots \}$ in Step {\it 3} and in {\it (b)} of Step {\it 4} in Algorithm \ref{algo:WT}.
\end{remark}

\begin{remark}
The output of Algorithm \ref{algo:TP} is not uniformly distributed over all elements of Graver basis.
The distribution depends on how to implement the randomness in Step {\it 3} and in {\it (b)} of Step {\it 4} in Algorithm \ref{algo:WT}.
\end{remark}

Algorithm \ref{algo:TP} allows us to uniformly sample graphs with the common degree sequence via Metropolis-Hastings algorithm with the Graver basis, 
with or without the restriction that graphs are simple.
It is done by constructing a connected Markov chain of graphs 
with the common degree sequence.
In each iteration, a primitive walk is randomly generated 
by Algorithm \ref{algo:TP}.
If the primitive walk is applicable, a new sample graph 
with the same degree sequence is obtained by adding the primitive walk, 
otherwise the primitive walk is rejected.
Note that Metropolis-Hastings algorithm does not require the uniformity
of the distribution of generated primitive walks.
As long as there is a positive probability of generating every element of the Graver basis,
the Metropolis-Hastings algorithm realizes uniform sampling of graphs with the common degree sequence.

\section{Numerical experiments}
\label{sec:examples}

In this section we present numerical experiments with elements of the Graver basis computed by Algorithm \ref{algo:TP} in Section \ref{sec:algorithm}.
The implementation of Metropolis-Hastings algorithm 
with Algorithm \ref{algo:TP} is done by Java 1.6.0
on Windows OS with Intel(R) Core(TM) i7-2829QM CPU@2.30GHz.

\subsection{A simulation with a small graph}
\label{subsec:example}
\begin{figure}[!h]
\vspace{-10mm}
\begin{minipage}{0.5\hsize}
\begin{center}
\includegraphics[]{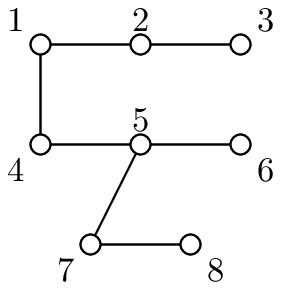}
\figcaption{Small graph $H_0$.}
\label{fig:mini_ex}
\end{center}
\end{minipage}
\begin{minipage}{0.5\hsize}
\begin{center}
\includegraphics[width=5cm]{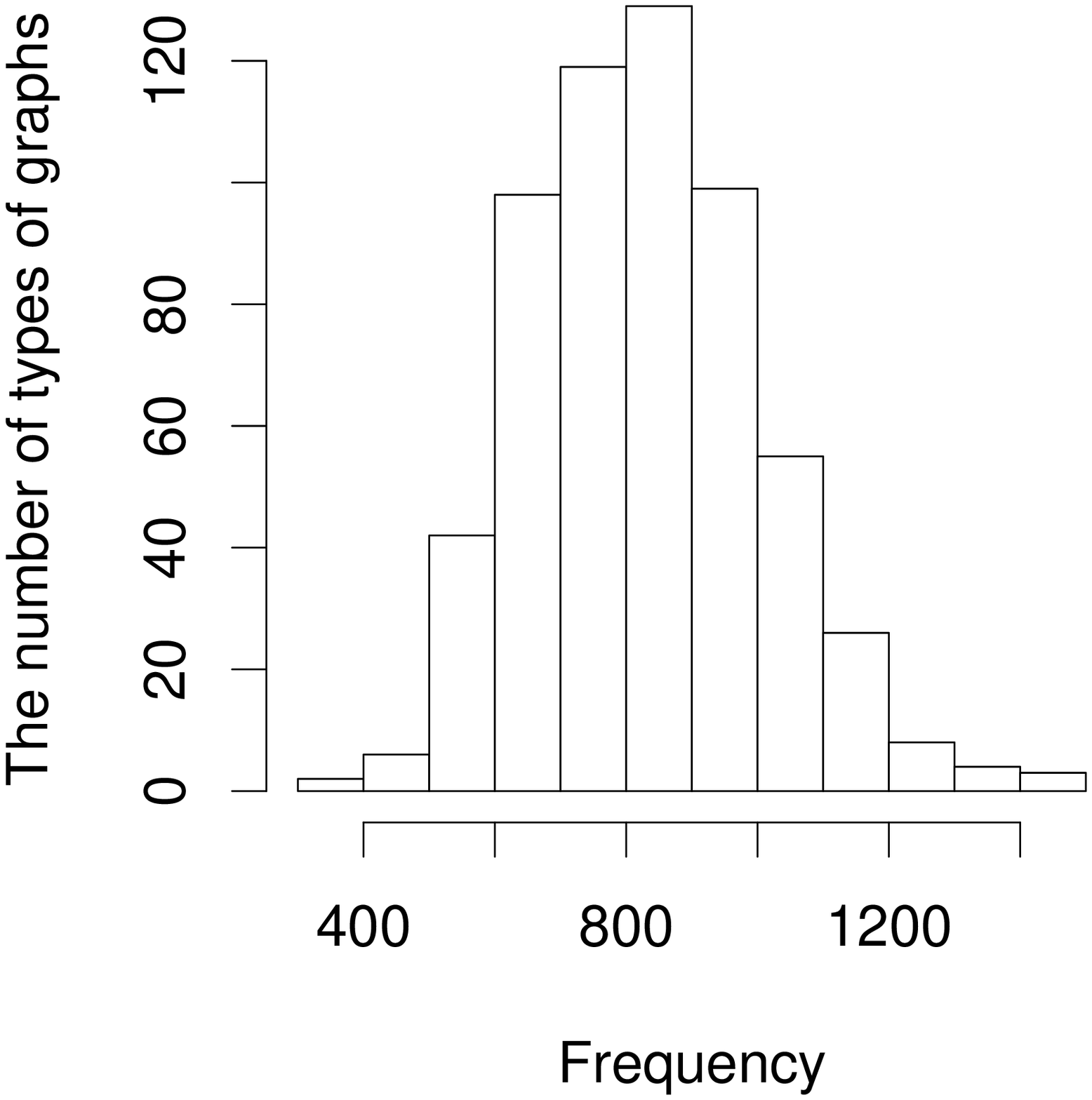}
\vspace{-3mm}
\figcaption{Histogram from sampling.}
\label{fig:hist}
\end{center}
\end{minipage}
\end{figure}

We run a Markov chain over the 
fiber containing a small graph $H_0$ in \figref{fig:mini_ex}.
The underlying graph $G=K_8$ is assumed to be complete with eight vertices.
By the Markov chain we sampled 510,000 graphs in the fiber, including 10,000 burn-in steps.
The number of types of obtained graphs in our chain is 591.
By enumeration we checked that 591 is actually the number of the elements of the fiber of $H_0$.
The histogram of this experiment is shown in Figure \ref{fig:hist}.
The horizontal axis expresses the frequency of each type of graph and the vertical axis expresses the number of types.
The mean of the number of appearances of each type is 829 and the standard deviation is 179.
This experiment shows that the algorithm samples each element of the fiber almost uniformly.

\subsection{The beta model for the food web data}
\label{subsec:real_data}

We apply Algorithm \ref{algo:TP} for testing of the real data, the
observed food web of 36 types of organisms in the Chesapeake Bay during
the summer. 
This data is available online at \cite{food_web}.
\citet{blitzstein-diaconis} analyzed essentially the same data set.
\vspace{-1mm}
\begin{figure}[!h]
\begin{center}
\includegraphics[width=8.5cm]{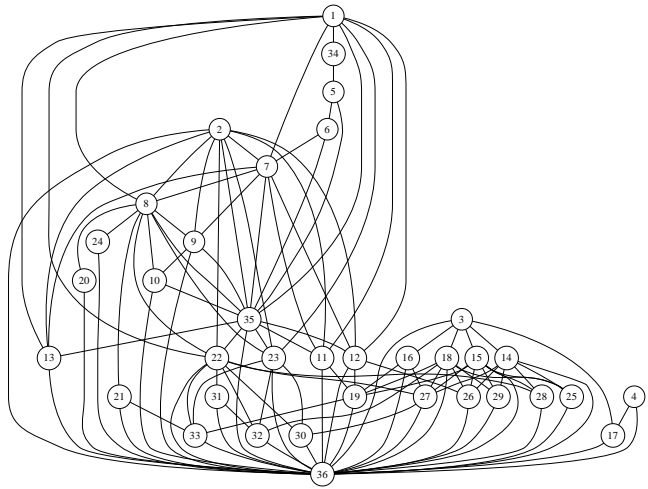}
\vspace{-1mm}
\figcaption{Food web for the Chesapeake Bay during the summer.}
\label{fig:food_web}
\end{center}
\end{figure}

The graph $H$ of the data is shown in \figref{fig:food_web}.
The vertices represent the types of organisms like blue
crab, bacteria etc., and the edges represent the relationship of one
preying upon the other. 
The degree sequence of $H$ is
\begin{eqnarray*}
&(9,10,6,2,3,3,9,11,6,4,6,7,5,7,8,4,3,8,\\
& \hspace{4mm}7,2,3,11,8,2,4,5,7,4,4,4,3,5,5,2,14,29). 
\end{eqnarray*}
Although there is a self loop at the vertex 19 in the observation, 
we ignored it for simplicity.

\begin{figure}[!h]
\vspace{-8mm}
\begin{center}
\includegraphics[width=5.5cm]{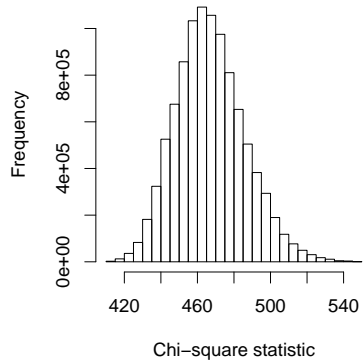}
\vspace{-4mm}
\figcaption{Histogram of chi-square statistic.}
\label{fig:hist_chisq}
\end{center}
\end{figure}

We set the beta model (\ref{beta-model}) in Section \ref{sec:beta-model} with $n_{ij}=1$ for each edge $\{i,j\}$ as the null hypothesis.
Then the probability of $H$ is described as
\begin{align}
 \label{beta_null}
 P(H) \propto \frac{\prod_{i \in V} \alpha_i^{d_i}}{\prod_{\{ i,j\} \in E} (1+\alpha_i \alpha_j )}.
\end{align}
Parameter $\alpha_i~(i \in V)$ is interpreted as the value of organism represented by the vertex $i$ as a food to other organisms.
Then the beta model (\ref{beta_null}) implies that a vertex $i$ with large $\alpha _i$ is likely to be connected to many edges.
Let $P \in (\ref{beta_null})$ mean that $P$ can be expressed by (\ref{beta_null}) for a set of parameters $\{ \alpha _i \} _{i \in V}$.
Consider now the statistical hypothesis testing problem
\[
H_0:P \in (\ref{beta_null}) ~\text{versus} ~
H_1: P \notin  (\ref{beta_null}).
\]
Starting from the graph in \figref{fig:food_web}, we construct 
a Markov chain of 10,100,000 graphs including 100,000 burn-in steps
 and compute the chi-square statistic 
of each graph as a test statistic.
The running time of the calculation is 5 minutes and 4.8 seconds.
Using the maximum likelihood estimator, the chi-square value of observed graph $H$ is 477 and the histogram of the estimated distribution of the chi-square values is shown in \figref{fig:hist_chisq}. The approximate $p$-value is 0.286. 
This value is not so small and there is no evidence against the beta model (\ref{beta_null}).

\begin{figure}[!h]
\vspace{-8mm}
\begin{minipage}{0.5\hsize}
\begin{center}
\includegraphics[width=5.0cm]{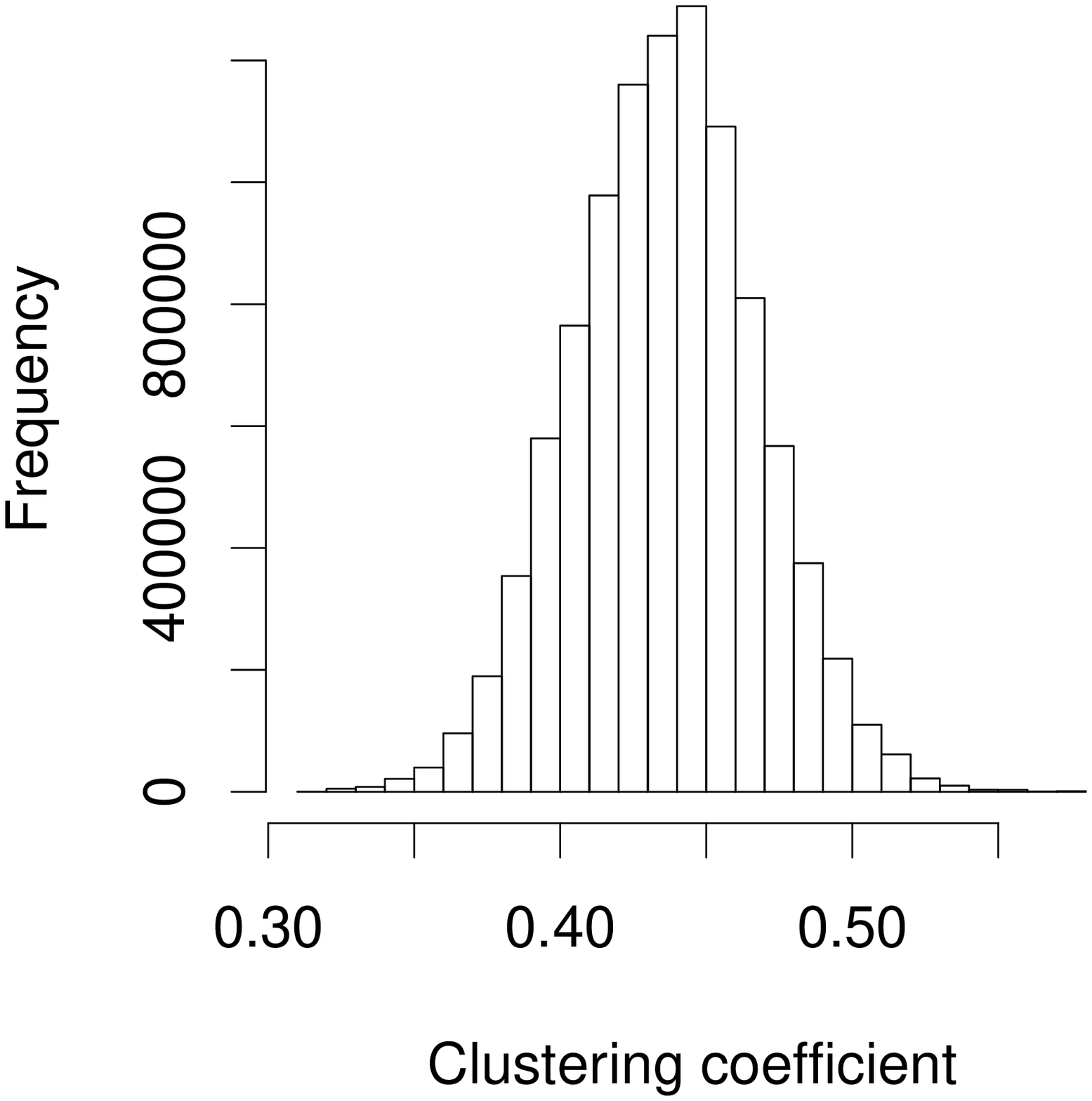}
\figcaption{Histogram of clustering coefficient.}
\label{fig:clco}
\end{center}
\end{minipage}
\begin{minipage}{0.5\hsize}
\begin{center}
\includegraphics[width=5.0cm]{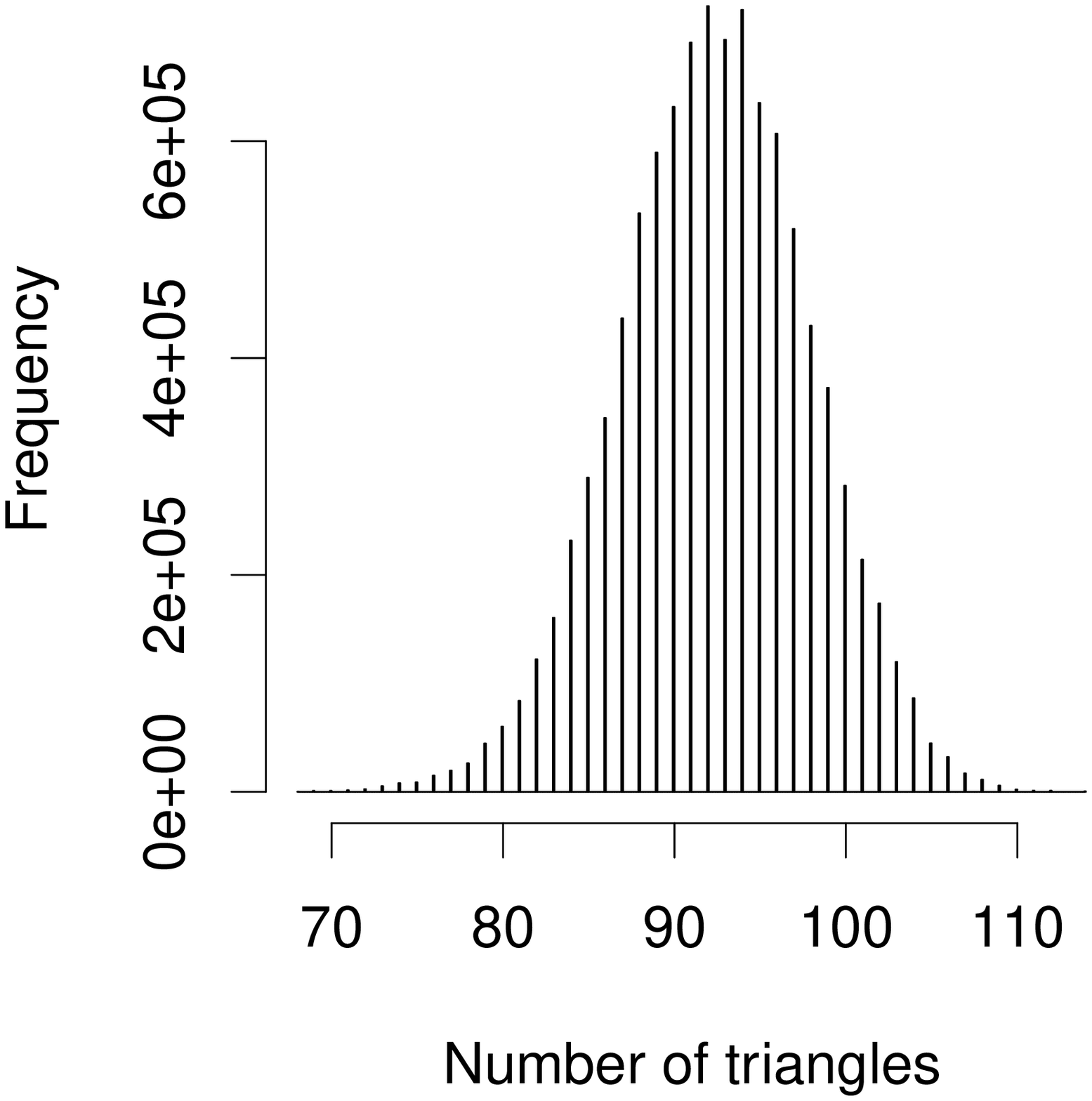}
\vspace{-3mm}
\figcaption{Histogram of number of triangles.}
\label{fig:ntri}
\end{center}
\end{minipage}
\end{figure}

Next we consider some other characteristics 
of the observed graph $H$ and graphs
obtained by the above Markov chain.
We compute their clustering coefficient
defined by \citet{watts-strogatz-1998}
and also count the number of triangles ($3$-cycles).
For the observed graph $H$, 
the values of clustering coefficient and the number of triangles
are 0.447 and 101, respectively.
For the sampled graphs,
the histograms are obtained as in \figref{fig:clco} and \ref{fig:ntri}
and their mean values are 0.436 and 92.4, respectively.
The differences between the actual values and the means of sampled graphs
are not large.
It suggests that these statistics agree with the beta model (\ref{beta_null}).

As mentioned in Section \ref{sec:introduction}
there are computer algebra systems such as 4ti2 (\citet{4ti2})
to compute the Graver basis.
However the whole Graver basis is huge and difficult to compute even for 
a moderate-sized graph like the real data above.
Algorithm \ref{algo:TP}, our adaptive algorithm, enables us
to perform the Markov chain Monte Carlo method for 
such a moderate-sized graph.

\section{Concluding remarks}
\label{sec:remarks}
In this paper we obtained a simple characterization of the Graver basis
for toric ideals arising from undirected graphs.  This Graver basis
allows us to perform the conditional test of the beta model for arbitrary
underlying graph.  Our characterization allows us to construct an
algorithm for sampling elements of the Graver basis, which is sufficient
for performing the conditional test.

By numerical experiments we confirmed that our procedure works well in practice.
We should mention that 
the sequential importance sampling method of \citet{blitzstein-diaconis} 
may work faster for the case of complete underlying graph.

If we allow multiple edges, then we do not need the Graver basis.  A minimal
Markov basis, which is often much smaller than the Graver basis, is sufficient
for connectivity of Markov chains.  Properties of Markov basis for the
$p_1$-model have been given in \citet{petrovic-etal-Urbana2010}.
It is of interest to study properties of minimal Markov bases for undirected graphs,
including the case of allowing self loops.


%
%

\bigskip
\noindent {\bf Acknowledgements}\quad 
We are very grateful to Hidefumi Ohsugi for valuable discussions.
We also thank two referees for their valuable and constructive comments.


\bibliographystyle{spbasic} 
\bibliography{graph-graver}

\appendix
\section{Proofs of Lemmas in Section \ref{subsec:characterization}}

\subsection{Proof of Lemma \ref{lemma:nc2}}
\begin{figure}[!h]
\begin{center}
\includegraphics[width=0.6\textwidth]{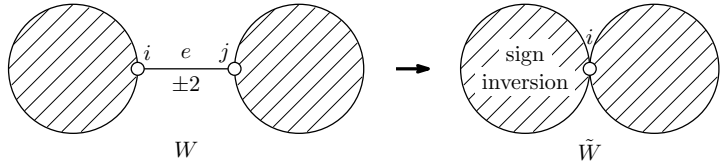}
\vspace{-2mm}
\figcaption{Contraction of an edge $e$.}
\label{fig:reductionforlemma}
\end{center}
\end{figure}

The contraction of the edge with its weight $\pm 2$ on $W$ is possible from Lemma \ref{lemma:nc1}.
We denote this edge by $e= \{ i, j\}$ as shown in \figref{fig:reductionforlemma}.
Suppose $\tilde{w}$ is not primitive.
Then there exists an proper subwalk $\tilde{w}^{\prime}$ of $\tilde{w}$.
If $i \notin V(\tilde{w}^{\prime})$,
$\tilde{w}^{\prime}$ is also a proper subwalk of $w$, a contradiction to the primitiveness of $w$.
Then $i \in V(\tilde{w}^{\prime})$.
However, a proper subwalk of $w$ is constructed by embedding $e$ into $\tilde{W}^{\prime}$.
Therefore, $\tilde{w}$ is primitive.
\qed

\subsection{Proof of Lemma \ref{lemma:nc3}}
\begin{figure}[!h]
\begin{center}
\includegraphics[width=0.7\textwidth]{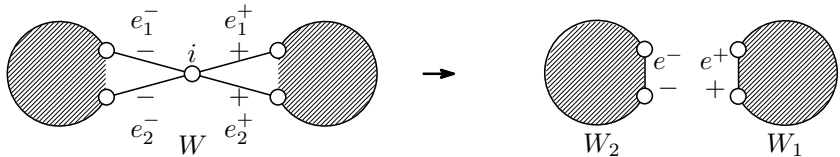}
\vspace{-2mm}
\figcaption{Separation of a vertex $i$.}
\label{fig:sepforlem}
\end{center}
\end{figure}

We consider the case that both positive and negative sides of $i$ correspond to (a) in \figref{fig:separation1}
and relevant edges are labeled as shown in \figref{fig:sepforlem}.
Suppose $w_1$ is neither primitive nor of length two.
Then  there exists a proper subwalk $w_1^{\prime}$ of $w_1$ on $W_1$.
If $e^+ \notin E(w_1^{\prime})$,
$w_1^{\prime}$ is also a proper subwalk of $w$, a contradiction to the primitiveness of $w$.
Then $e^+ \in E(w_1^{\prime})$.
Now $w_1^{\prime}$ is expressed as follows:
\begin{eqnarray*}
	w_1^{\prime}=(e_{i_1}, e_{i_2},\ldots, e_{i_k},e^+,e_{i_{k+1}},\ldots ,e_{i_s}).
\end{eqnarray*}
Then an even closed walk on $W$
\begin{eqnarray*}
	(e_{i_1}, e_{i_2},\ldots, e_{i_k},e_1^+,e_1^-,\ldots,e_2^-,e_2^+,e_{i_{k+1}},\ldots ,e_{i_s})
\end{eqnarray*}
is a proper subwalk of $w$.
This contradicts the primitiveness of $w$.
Therefore $w_1$ is primitive or of length two.
The cases of (b) and (c) in \figref{fig:separation1} are shown in the same way.
Note that it is easy to confirm the possibility of contraction after the step 1 in the case (b) from Lemma \ref{lemma:nc1} and then the primitiveness is guaranteed by Lemma \ref{lemma:nc2}.
By the same argument, the case of $w_2$ is confirmed.
\qed

\subsection{Proof of Lemma \ref{lemma:sc2}}
\begin{figure}[!h]
\begin{center}
\includegraphics[width=0.6\textwidth]{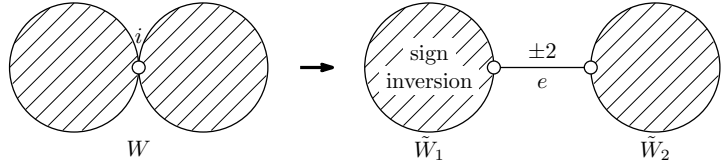}
\vspace{-2mm}
\figcaption{Insertion to a vertex $i$.}
\label{fig:redforlem}
\end{center}
\end{figure}

Let $e$ be the new edge appearing through the insertion to $i$ as shown in \figref{fig:redforlem}.
Suppose $\tilde{w}$ is not primitive.
Then there exists a proper subwalk $\tilde{w}^{\prime}$ of $\tilde{w}$.
If $e \notin E(\tilde{w}^{\prime})$,
$\tilde{w}^{\prime}$ is contained in $\tilde{W}_1$ or $\tilde{W}_2$.
Then $\tilde{w}^{\prime}$  or its reverse
becomes a proper subwalk of $w$.
This contradicts  the primitiveness of $w$.
Hence $e \in E(\tilde{w}^{\prime})$.
Then we can construct a proper subwalk of $w$ by removing $e$ from $\tilde{w}^{\prime}$ and reversing the weights of edges belonging to $E(w_1)$, a contradiction to the primitiveness of $w$.
Therefore, $\tilde{w}$ is primitive.
\qed

\subsection{Proof of Lemma \ref{lemma:sc1}}
\begin{figure}[!h]
\begin{center}
\includegraphics[width=0.7\textwidth]{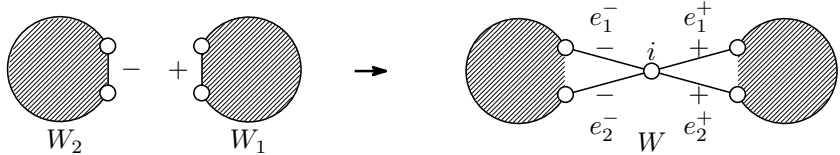}
\vspace{-2mm}
\figcaption{Binding of $W_1$ and $W_2$.}
\label{fig:combine}
\end{center}
\end{figure}

Let $i$ be the new vertex appearing through the binding.
We consider the case that both positive and negative sides of $i$  correspond to (a) in \figref{fig:separation1} and relevant edges are labeled as shown in \figref{fig:combine}.
Other cases are shown in the same way.
Suppose $w$ is not primitive.
Then there exists a proper subwalk $w^{\prime}$ of $w$.
Here we choose a primitive walk as $w^{\prime}$.
If $i \notin V(w^{\prime})$, $w^{\prime}$ is also a proper subwalk of $w_1$ or $w_2$.
Then $i \in V(w^{\prime})$.
This implies that all four edges connected to $i$ appear in $w^{\prime}$.
Let us consider the separation of $i$ to $W^{\prime}$.
Then the resulting two weighted graphs $W_1^{\prime},W_2^{\prime}$ are primitive from Lemma \ref{lemma:nc3}.
Furthermore at least one of $w_i^{\prime}\ (i=1,2)$ is a proper subwalk of $w_i$, a contradiction to the primitiveness of $w_i$.
Therefore $w$ is primitive.
\qed

%
%
%
%
%

\end{document}